\documentclass[12pt]{amsart}
\usepackage{amssymb, amstext, amscd, amsmath}
\usepackage{fullpage}

\theoremstyle{plain}

\theoremstyle{definition}

\theoremstyle{remark}

\begin{document}

\title{Cocycle perturbation on Banach algebra}
\title{Cocycle perturbation on Banach algebra}

\author[Y. J. Wu]{YU Jing Wu}

\address{Tianjin Vocational Institute \\Tianjin 300160\\P.R. CHINA}

\email{wuyujing111@yahoo.cn}

\author[L. Y. Shi]{Luo Yi Shi}

\address{Department of Mathematics\\Tianjin Polytechnic University\\Tianjin 300160\\P.R. CHINA}

\email{sluoyi@yahoo.cn}

%
%

\thanks{Supported by NCET(040296), NNSF of China(10971079)}

\date{}

%
%



\begin{abstract}

Let $\alpha$ be a flow on a Banach algebra $\mathfrak{B}$,
and $t\longmapsto u_t$ a continuous function on $\mathbb{R}$ into the group of invertible
elements of $\mathfrak{B}$ such that $u_s\alpha_s(u_t )=u_{s+t}, s, t \in \mathbb{R}$. Then
$\beta_t=$Ad$u_t\circ\alpha_t, t\in \mathbb{R}$ is also a flow on $\mathfrak{B}$.
$\beta$ is said to be a cocycle perturbation of $\alpha$.
 We show that if $\alpha,\beta$ are two flows on nest algebra (or quasi-triangular algebra), then
$\beta$ is a cocycle perturbation of $\alpha$.  And the flows on
nest algebra (or quasi-triangular algebra) are all uniformly
continuous. \end{abstract}

\maketitle

\section{Introduction}

In the quantum mechanics of particle systems with an infinite number
of degrees of freedom, an important problem is to study the
differential equation
$$\frac{d\alpha_t(A)}{dt}=S\alpha_t(A)$$
under variety of circumstances and assumptions. In each instance the
$A$ corresponds to an observable, or state, of the system and is
represented by an element of some suitable Banach algebra
$\mathfrak{B}$. And $S$ is an operator on $\mathfrak{B}$,
$\{\alpha_t\}_{t\in\mathbb{R}}$ is a group of bounded automorphisms
on $\mathfrak{B}$. The Function
$$t\in \mathbb{R}\longmapsto \alpha_t(A)\in\mathfrak{B} $$
describes the motion of $A$. The dynamics are given by solutions of the differential equation which
respect certain supplementary conditions of continuity. Thus it is worth to
study the group of bounded automorphisms on
$\mathfrak{B}$. More details see [1].

A {\it flow} $\alpha$ on $\mathfrak{B}$ is a group homomorphism of the real
line $ \mathbb{R}$ into the group of bounded automorphisms on
$\mathfrak{B}$ (i.e. $t\longmapsto\alpha_t$) such that
 $$\lim _{t\rightarrow t_0}||\alpha_t(B)-\alpha_{t_0}(B)||=0$$
  for
each $t_0\in\mathbb{R}$ and each $B\in\mathfrak{B}$. If there exists a $h\in \mathfrak{B}$ such that
$\alpha_t(B)=e^{th}Be^{th},\forall B\in\mathfrak{B},t\in\mathbb{R}$, then we call $\alpha$  an {\it inner flow}.
 We say that a flow $\alpha$ is {\it uniformly continuous} if
$$\lim _{t\rightarrow t_0}||\alpha_t-\alpha_{t_0}||=0$$
 for
each $t_0\in\mathbb{R}$.

If $\alpha$ is a flow on $\mathfrak{B}$ and if $u$ is a continuous map of $ \mathbb{R}$
into the group of invertible
elements $G(\mathfrak{B})$ of $\mathfrak{B}$ such that
$u_s\alpha_s(u_t )=u_{s+t}, s, t \in \mathbb{R}$, then we call $u=(u_t)_{t\in\mathbb{R}}$
 an {\it $\alpha$-cocycle} for $(\mathfrak{B},\mathbb{R},\alpha)$.
 Let
$\beta_t=$Ad$u_t\circ\alpha_t, t\in \mathbb{R}$ (i.e. $\beta_t(B)=u_t\alpha_t(B)u_t^{-1}, \forall B\in\mathfrak{B}$),
then $\beta$
  is also a flow on $\mathfrak{B}$,
and $\beta$ is said to be a {\it cocycle perturbation }of $\alpha$.

If $\alpha$ is a flow on $\mathfrak{B}$, let
$D(\delta_\alpha)$ be composed of those $B\in\mathfrak{B}$ for which there exists a $A\in\mathfrak{B}$ with the
property that
$$A=\lim _{t\rightarrow 0}\frac{\beta_t(B)-B}{t}.$$
Then $\delta_\alpha$ is a linear operator on $D(\delta_\alpha)$ defined by
$\delta_\alpha(B)=A$. We call $\delta_\alpha$ the {\it infinitesimal generator} of
$\alpha$. By [1,Proposition 3.1.6],
 $\delta_\alpha$ is a closed derivation. i.e. the domain $D(\delta_\alpha)$ is a dense
subalgebra of $\mathfrak{B}$ and $\delta_\alpha$ is closed as a linear operator
on $D(\delta_\alpha)$ and satisfies
$\delta_\alpha(AB)=\delta_\alpha(A)B+A\delta_\alpha(B),$
for $A,B\in D(\delta_\alpha)$. We call $\beta$
an {\it inner perturbation} of $\alpha$, if $\alpha,\beta$ are two flows on $\mathfrak{B}$,
 $D(\delta_\alpha)=D(\delta_\beta)$ and there exists  $h\in\mathfrak{B}$ such that
 $\delta_\beta=\delta_\alpha+$ad$ih$ (where $i$ is the imaginary unit, and ad$ih(B)\triangleq i(hB-Bh), \forall B\in\mathfrak{B}$).
  Moreover $D(\delta_\alpha)=\mathfrak{B}$ if and only if $\alpha$ is uniformly continuous.
More details see [1,2].

The problem we considered is classifying cocycle of flows on Banach algebras. Such a problem
has been considered in the $C^*$algebra cases, notably by Kishimoto[3-10]. We refer the reader to [3] for a detailed
study of the general results concerning cocycles and invariants for cocycle perturbation
 of flows on $C^*$-algebras. In particular, it is shown there that
let $\alpha$ be a flow on a $C^*$-algebra $\mathfrak{A}$, if
$u=(u_t)_{t\in\mathbb{R}}$  a $\alpha$-cocycle for $(\mathfrak{A},\mathbb{R},\alpha)$, and $u$ is differentiable (i.e.
$\lim_{t\rightarrow t_0}\frac{u_t-u_{t_0}}{t-t_0}$ exist for any $t_0\in \mathbb{R}$) and $h=-i(du_t/dt)|_{t=0}\in \mathfrak{A}$,
then the infinitesimal generator of the flow
 $\beta_t= Adu_t\circ\alpha_t$ is given by $\delta_\beta=\delta_\alpha+$ad$ih$. i.e. $\beta$ is an inner
 perturbation of $\alpha$. Moreover for any $\alpha$-cocycle $u=(u_t)_{t\in\mathbb{R}}$,
there is
  $w\in G(\mathfrak{A})$ and a differentiable $\alpha$-cocycle $v=(v_t)_{t\in\mathbb{R}}$
 (i.e. $v_t$ is an $\alpha$-cocycle and differentiable)
  such that $u_t=wv_t\alpha_t(w^{-1})$.
In section 2, we consider the cocycle of flows on Banach algebras and obtain some similar results to [3].

 It is well-known that
 a flow $\alpha$ on $\mathfrak{B}$ may not be uniformly continuous even
 if $\mathfrak{B}$ is a $C^*$algebra [1,2].
 In section 3, we study the flows on nest algebra $\tau(\mathcal{N})$ and
quasi-triangular algebra $Q\tau(\mathcal{N})=\tau(\mathcal{N})+K$
[11]. We recall that a {\it nest} $\mathcal{N}$ is a chain of closed
subspaces of a Hilbert space $\mathfrak{H}$ containing $\{0\}$ and
$\mathfrak{H}$ which is in addition closed under taking arbitrary
intersections and closed spans. The {\it nest algebra}
$\mathcal{T(N)}$ associated with $\mathcal{N}$
 is the set of all $T\in B(\mathfrak{H})$ which leave each element of the nest invariant.
For instance, if $\mathfrak{H}$ is separable with orthonormal basis $\{e_n\}_{n=1}^\infty$ and
$\mathfrak{H}_n$=span$\{e_1, \cdots, e_n\}$, then
$\mathcal{N}=\{0, \mathfrak{H}\}\cup\{\mathfrak{H}_n\}_{n=1}^\infty$ is a nest.
In this case, $\mathcal{T(N)}$ is simply the set of all operators whose matrix representation with
respect to this basis is upper triangular.
It is obvious,
 $\tau(\mathcal{N})$ and
$Q\tau(\mathcal{N})$ are typical  Banach algebras.
 We obtain that
all of the flows on $\tau(\mathcal{N})$ (or $Q\tau(\mathcal{N})$)
are uniformly continuous. Moreover all of the flows are cocycle perturbation
to each other.

\vskip1cm
\section{Cocycle perturbations}

Let $\mathfrak{B}$ be a Banach algebra, $\alpha$ be a flow on
$\mathfrak{B}$ and $u$ an $\alpha$-cocycle, then
$\beta_t=$Ad$u_t\circ\alpha_t, t\in \mathbb{R}$
is a cocycle perturbation of $\alpha$. In this section, we obtain
that $\beta$ is an inner perturbation of $\alpha$ if and only if $u$
is differentiable (Theorem 2.4.). Moreover for any $\alpha$-cocycle $u$ ,
there is a differentiable $\alpha$-cocycle $v$
and a inventible element $w$ in
$\mathfrak{B}$ such that $u_t=wv_t\alpha_t(w^{-1})$ (Theorem 2.8.).

We will use the following Lemmas [1, Proposition 3.1.3, Theorem 3.1.33].

{\bf Lemma 2.1$^{[1]}$:} Let $\{\alpha_t\}_{t\in\mathbb{R}}$ be a flow
on the Banach algebra
$\mathfrak{B}$. Then there exists an $M \geq 1$ and $\xi\geq$inf$_{t\neq
0}$($t^{-1}$log$||\alpha_t||$) such that
$||\alpha_t||\leq Me^{\xi |t|}.$

{\bf Lemma 2.2$^{[1]}$:} Let $\alpha$ be a flow on a Banach algebra $\mathfrak{B}$ with infinitesimal
 generator $\delta_\alpha$.  For each $P\in\mathfrak{B}$ define the bounded
derivation $\delta_P$ by $D(\delta_P)=\mathfrak{B}$ and
$\delta_P(B)=i[P,B]\triangleq i(PB-BP), \forall B\in\mathfrak{B}$. Then
$\delta+\delta_P$ generates a flow on $\mathfrak{B}$ given by
$$\alpha_t^P(B)=\alpha_t(B)+\sum_{n\geq 1}i^n\int_0^t dt_1\int_0^{t_1}
dt_{2}\cdots\int_0^{t_{n-1}}dt_n
[\alpha_{t_n}(P),[\cdots[\alpha_{t_1}(P),\alpha_t(B)]]], \forall B\in\mathfrak{B}, t\in\mathbb{R}.$$

{\bf Lemma 2.3:} Let $\mathfrak{B}$ be a Banach algebra with unit $\mathbf{1}$,  $\alpha$ be a flow on
$\mathfrak{B}$ and let $\delta$ denote the infinitesimal generator of
$\alpha$. Furthermore, for each $P\in\mathfrak{B}$ define  $\delta_P$ as in Lemma 2.2.

Then $\delta+\delta_P$ generates a flow on $\mathfrak{B}$ given by
$$\alpha_t^P(B)=\alpha_t(B)+\sum_{n\geq 1}i^n\int_0^t dt_1\int_0^{t_1}
dt_{2}\cdots\int_0^{t_{n-1}}dt_n
[\alpha_{t_n}(P),[\cdots[\alpha_{t_1}(P),\alpha_t(B)]]].$$
Moreover, one has $$\alpha_t^P(B)=u^P_t\alpha_t(B){(u^P_t})^{-1},$$
 where $u_t^p$ is a
one-parameter family of invertible elements, determined by
$$u^P_t=1+\sum_{n\geq 1}i^n\int_0^t dt_1\int_0^{t_1}
dt_{2}\cdots\int_0^{t_{n-1}}dt_n\alpha_{t_n}(P)\cdots\alpha_{t_1}(P)$$

which satisfies the $\alpha$-cocycle relation
$$u^P_{t+s}=u^P_t\alpha_t(u^P_s).$$
All integrals converge in the strong topology.  The integrals define
norm-convergent series of bounded operators and there exists an $M \geq 1$, $\xi\geq$inf$_{t\neq
0}$($t^{-1}$log$||\alpha_t||$) such that

$||\alpha_t^P(B)-\alpha_t(B)||\leq Me^{\xi |t|}(e^{M|t|||P||} -
1),||u^P_t-\mathbf{1}||\leq Me^{\xi |t|}(e^{M|t|||P||} - 1)$ .

{\bf Proof :}   The first statement of the proposition can be
obtained from Lemma 2.2. We just give the proof of the last statement of the Lemma.

 There exists an $M \geq 1$ and $\xi\geq$inf$_{t\neq
0}$($t^{-1}$log$||\alpha_t||$) such that
$||\alpha_t||\leq Me^{\xi |t|},$ by Lemma 2.1. Let
$$M_t=\left\{
\begin{array}{ll}
   Me^{\xi |t|} &\mbox{$|t|>1$;}\\
   M   &\mbox{$|t|\leq 1$.}
\end{array}
\right.$$
Next we consider $u^P_t$ defined by the series.
The $n$-th term in this series is well defined and has norm less than
$\frac{|t|^n}{n!}M_t^n||P||^n$.  Thus $u^P_t$ is a norm-continuous
one-parameter family of elements of $\mathfrak{A}$ with $u^P_0=1$
and $||u_t||\leq e^{|t|M_t||P||}$. Consequently,  $u^P_t$  is
invertible for all $t\in[-t_0,t_0]$ for some $t_0>0$ and
$(u^P_t)^{-1}$ is a norm-continuous one-parameter family of elements
of $\mathfrak{B}$ for all $t\in[-t_0,t_0]$.

 Next one has
$$\frac{du^P_t}{dt}=iu^P_t\alpha_t(P),$$
$$\lim_{t\rightarrow 0}(\frac{u^P_t-1}{t})=\frac{du^P_t}{dt}|_{t=0}=iP.$$
Hence,
$$\lim_{t\rightarrow 0}(\frac{(u^P_t)^{-1}-1}{t})=
\lim_{t\rightarrow 0}(u^P_t)^{-1}\frac{1-u^P_t}{t}=-iP.$$
To establish the $\alpha$-cocycle relation we first note that
$$\frac{du^P_{t+s}}{ds}=iu^P_{t+s}\alpha_{t+s}(P)$$
and $u^P_{t+s}|_{s=0}=u^P_t$.

But $\alpha_t(u^P_s)=u^{\alpha_t(P)}_s$
and hence
$$\frac{d}{ds}u^P_t\alpha_t(u^P_s)=
iu^P_tu^{\alpha_t(P)}_s\alpha_s(\alpha_t(P))=
iu^P_t\alpha_t(u^P_s)\alpha_{t+s}(P).$$
Moreover, $u^P_t\alpha_t(u^P_s)|_{s=0}=iu^P_t\alpha_t(p)$. Thus $s\longmapsto
u^P_{t+s}$ and $s\longmapsto u^P_t\alpha_t(u^P_s) $ satisfy the same
first-order differential equation and boundary condition for each
$t\in\mathbb{R}$. Therefore, the two functions are equal and can be
obtained by iteration of the integral equation
$$X_t(s)=u^P_t+i\int_0^sds'X_t(s')\alpha_{t+s'}(P).$$
Hence, $u^P_{t+s}=u^P_t\alpha_t(u^P_s),t,s\in\mathbb{R}$. Since $u^P_t$
is invertible for all $t\in[-t_0,t_0]$ for some $t_0>0$, we obtain that
$u^P_t$ is a norm-continuous one-parameter family of inventible
elements.
Thus $t\mapsto u^P_t\alpha_t(B)(u^P_t)^{-1}$ defines a flow $\beta_t$ on
$\mathfrak{B}$.

Let $\tilde{\delta}$ denote the infinitesimal generator of $\beta$.
Next we will prove $\tilde{\delta}=\delta+\delta_P$.

Choose $A\in D(\delta+\delta_P)$, one has
\begin{equation*}
\delta(A)=\lim_{t\rightarrow 0}(\frac{\alpha_t(A)-A}{t}),
\end{equation*}
\begin{eqnarray*}
\tilde{\delta}(A)&=&\lim_{t\rightarrow 0}(\frac{\beta_t(A)-A}{t})\\
&=&\lim_{t\rightarrow 0}(\frac{u^P_t\alpha_t(A)(u^P_t)^{-1}-u^P_tA(u^P_t)^{-1}}{t}\\
&+&\frac{u^P_tA(u^P_t)^{-1}-u^P_tA}{t}+\frac{u^P_tA-A}{t})\\
&=&(\delta+\delta_P)(A).
\end{eqnarray*}

Similarly, if $A\in D(\tilde{\delta})$, we
obtain$\tilde{\delta}(A)=(\delta+\delta_P)(A)$. And then
$\tilde{\delta}=\delta+\delta_P$.

Thus one must have
$\alpha_t^P(B)=\beta_t(B)=u^P_t\alpha_t(B){(u^P_t})^{-1},$ for any $B\in\mathfrak{B}$
by [1, Theorem 3.1.26].
Finally the estimates on $\alpha_t^P(B)-\alpha_t(B)$ and
$u^P_t-\mathbf{1}$ are
straightforward.

{\bf Theorem 2.4 :} Let $\alpha$ be a flow on $\mathfrak{B}$, $(u_t)_{t\in\mathbb{R}}$ be an $\alpha$-cocycle
, $\beta_t=$Ad$u_t\circ\alpha_t$. Then $\beta$ is an  inner perturbation of $\alpha$ if and only if $u_t$ is
is differentiable.

{\bf Proof :}
$\Rightarrow$ It follows immediately from Lemma 2.3.

$\Leftarrow$ If $u_t$ is an $\alpha$-cocycle and differentiable  with
$h=-i(du_t/dt)|_{t=0}\in \mathfrak{B}$, then $u_t$ is given by
$$u_t=1+\sum_{n\geq 1}i^n\int_0^t dt_1\int_0^{t_1}
dt_{2}\cdots\int_0^{t_{n-1}}dt_n\alpha_{t_n}(P)\cdots\alpha_{t_1}(P).$$
Hence, $\beta$ is an  inner perturbation of $\alpha$ by lemma 2.3.

{\bf Corollary 2.5 :} Adopt the assumptions of the Lemma 2.3  and also
assume that $\alpha_t$ is an inner flow. i.e. there exists $h\in\mathfrak{B}$
such that $\alpha_t(A)=e^{ith}Ae^{-ith}, \forall A\in\mathfrak{B},t\in\mathbb{R}$.

Then
$$\alpha^P_t(A)=\Gamma_t^PA(\Gamma_t^P)^{-1}, u_t^P=\Gamma_t^Pe^{-ith}, $$
where $u_t^P$ is defined as in Lemma 2.3 and $\Gamma_t^P=e^{it(h+P)}$. i.e. $\alpha^P_t$ is an inner flow.

{\bf Proof :}
If $\Gamma_t^P=e^{it(h+P)}$ and $X_t=\Gamma_t^Pe^{-ith}$, then
$$\frac{dX_t}{dt}=i\Gamma_t^PPe^{-ith}=iX_t\alpha_t(P)$$
and $X_0=\mathbf{1}$. Thus, $X_t$ is the unique solution of the
integral equation
$$X_t=\mathbf{1}+i\int_0^tdsX_s\alpha_s(P).$$
This solution can be obtained by iteration and one finds $X_t=u_t^P$, where
$u_t^P$ is defined as in Lemma 2.3. And$\alpha_t^P(A)=u^P_t\alpha_t(A){(u^P_t})^{-1}=
\Gamma_t^PA(\Gamma_t^P)^{-1}$.

 Next we obtain that every $\alpha$-cocycle  is similar to a
differentiable $\alpha$-cocycle.

{\bf Definition 2.6:} Let $\alpha$ be a flow on $\mathfrak{B}$. $A\in\mathfrak{B}$ is called an analytic element
for $\alpha$ if there exist a analytic function $f:~~\mathbb{C}\rightarrow \mathfrak{B}$ such that
$f(t)=\alpha_t(A)$ for $t\in\mathbb{R}$.

{\bf Lemma 2.7:} Let $\alpha$ be a flow on Banach algebra
$\mathfrak{B}$ and $M,\xi$ are constants such that $||\alpha_t||\leq
Me^{\xi |t|}$. For $A\in\mathfrak{B}$, define
$$A_n=\sqrt{\frac{n}{\pi}}\int\alpha_t(A)e^{-nt^2-\xi t}dt,~~~n=1,2,\cdots.$$
Then each $A_n$ is an entire analytic element for
$\alpha_t$, $||A_n||<M||A||$ for all n; and $A_n\rightarrow A$ in
the weak topology as $n\rightarrow\infty$. In particular, the
$\alpha$ analytic elements form a normal-dense subspace of
$\mathfrak{B}$.

{\bf Proof :} $$f_n(z)=\sqrt{\frac{n}{\pi}}\int\alpha_t(A)e^{-n(t-z)^2-\xi (t-z)}dt$$

is well defined for all $z\in\mathbb{C}$,
since $t\mapsto e^{-n(t-z)^2}\in\mathbf{L}^1(\mathbb{R})$
for each $z\in\mathbb{C}$.

For $z=s\in\mathbb{R}$, we have
\begin{eqnarray*}
f_n(s)&=&\sqrt{\frac{n}{\pi}}\int\alpha_t(A)e^{-n(t-s)^2-\xi (t-s)}dt\\
&=&\sqrt{\frac{n}{\pi}}\int\alpha_{t+s}(A)e^{-nt^2-\xi t}dt=\alpha_s(A_n).
\end{eqnarray*}
But for $\eta\in\mathfrak{B}^*$ we have
$$\eta(f_n(z))=\sqrt{\frac{n}{\pi}}\int\eta(\alpha_t(A))e^{-n(t-z)^2-\xi (t-z)}dt.$$
Since $|\eta(\alpha_t(A))|\leq M||\eta||||A||$, it follows from the Lebesgue dominated convergence
theorem that $t\mapsto\eta(f_n(z))$ is analytic. Hence each $A_n$ is analytic for $\alpha_t(A)$.

Next, one derives the estimate
$$||A_n||\leq M||A||\sqrt{\frac{n}{\pi}}\int e^{-n(t)^2} dt\leq M||A||.$$
Next note that
$$\int e^{-nt^2-\xi t}dt=e^{\frac{\xi^2}{4n^2}}\sqrt{\frac{\pi}{n}}.$$
Hence
$$\eta(A_n-A)=\sqrt{\frac{n}{\pi}}\int e^{-nt^2-\xi t}\eta(\alpha_t(A)-
e^{-\frac{\xi^2}{4n^2}}A)dt$$
 for all $\eta\in\mathfrak{B}^*$. But for any $\varepsilon>0$ we may choose $\delta>0$
 such that $t<|\delta|$ implies $|\eta(\alpha_t(A)-A)|<\varepsilon$. Further, we can choose $N$ large enough that
 $\frac{\xi}{2N}<\frac{\delta}{2}$.
 It follows that if $n>N$ we have

\begin{eqnarray*}
|\eta(A_n-A)|&\leq&\sqrt{\frac{n}{\pi}}\int_{|t|<\delta} e^{-nt^2-\xi t}|\eta(\alpha_t(A)-
e^{-\frac{\xi^2}{4n^2}}A)|dt\\
&+&
\sqrt{\frac{n}{\pi}}\int_{|t|\geq\delta} e^{-nt^2-\xi t}|\eta(\alpha_t(A)-
e^{-\frac{\xi^2}{4n^2}}A)|dt
\end{eqnarray*}
On the other hand,

\begin{eqnarray*}
\sqrt{\frac{n}{\pi}}\int_{|t|<\delta} e^{-nt^2-\xi t}|\eta(\alpha_t(A)-
e^{-\frac{\xi^2}{4n^2}}A)|dt
\end{eqnarray*}
\begin{eqnarray*}
&\leq&\sqrt{\frac{n}{\pi}}\int_{|t|<\delta} e^{-nt^2-\xi t}|\eta(\alpha_t(A)-
A)|dt\\
&+&\sqrt{\frac{n}{\pi}}\int_{|t|<\delta} e^{-nt^2-\xi t}|\eta(A)(1-e^{-\frac{\xi^2}{4n^2}})|dt\\
& < &\varepsilon e^{\frac{\xi^2}{4n^2}}
 +|1-e^{-\frac{\xi^2}{4n^2}}|\cdot ||\eta||\cdot ||A||Me^{\frac{\xi^2}{4n^2}}.\\
\end{eqnarray*}
And
\begin{eqnarray*}
\sqrt{\frac{n}{\pi}}\int_{|t|\geq\delta} e^{-nt^2-\xi t}|\eta(\alpha_t(A)-
e^{-\frac{\xi^2}{4n^2}}A)|dt
\end{eqnarray*}
\begin{eqnarray*}
&\leq&\sqrt{\frac{n}{\pi}}||\eta||\cdot ||A||M\int_{|t|\geq\delta}e^{-nt^2}dt\\
&+&\sqrt{\frac{n}{\pi}}e^{-\frac{\xi^2}{4n^2}}||\eta||\cdot ||A||\int_{|t|\geq\delta}
e^{-nt^2-\xi t}dt.
\end{eqnarray*}

So $|\eta(A_n-A)|\rightarrow 0$ when $n\rightarrow\infty$, for any $\eta\in\mathfrak{B}^*$.
Finally note that the norm closure and the weak closure of convex set are the same, hence
the
$\alpha$ analytic elements form a normal-dense subspace of
$\mathfrak{B}$.

{\bf Theorem 2.8:} If $u$ is an $\alpha$-cocycle for
$\mathfrak{B}$. Then, given $\varepsilon>0$, there is a differentiable
$\alpha$-cocvcle $v$ and a $w\in G(\mathfrak{B})$ such that
$$||w-\mathbf{1}||<\varepsilon, u_t=wv_t\alpha_t(w^{-1}).$$

{\bf Proof :} The proof is similar to [3,Lemma 1.1] and is omitted.

Given two flows $\alpha$ and $\beta$ on a unitary
Banach algebra $\mathfrak{B}$ we say that $\beta$ is a {\it conjugate} to
$\alpha$ if there is a bounded automorphism $\sigma$ of
$\mathfrak{B}$ such that $\beta=\sigma\alpha\sigma^{-1}$.
Conjugate, cocycle perturbation and inner perturbation define three
equivalence relations. We say that $\beta$ is {\it cocycle-conjugate} to
$\alpha$ if there is a bounded  automorphism $\sigma$ of $\mathfrak{B}$ such that
$\beta$ is a cocycle perturbation of $\sigma\alpha\sigma^{-1}$. This
also defines an equivalence relation among the flows.

We say that $\alpha$ is {\it approximately inner} if there is a sequence
$\{h_n\}$ in $ \mathfrak{B}$ such that $\alpha_t=$ limAd $e^{th_n}$,
i.e., $\alpha_t(A)= $$\lim$Ad$e^{th_n}(A)$ for every $t\in
\mathbb{R}$ and $A\in \mathfrak{B}$, or equivalently, uniformly
continuous in $t$ on every compact subset of $\mathbb{R}$ and every
$A\in\mathfrak{B}$.
 A flow on a Banach algebra $\mathfrak{B}$ is said to be
{\it asymptotically inner} if there is a continuous function $h$ of
$\mathbb{R}_+$
 into $\mathfrak{B}$ such that
 $\lim_{s\rightarrow \infty}\max_{|t\leq 1|}\parallel\alpha_t(A)-$Ad$e^{th(s)}(A)\parallel=0$
 for any $A\in\mathfrak{B}$.

{\bf Corollary 2.9:} Let  $\alpha$ and $\beta$  be two flows
on Banach algebra $\mathfrak{B}$. Then the following conditions
are equivalent:

(i). $\beta$ is cocycle-conjugate to $\alpha$.

 (ii). $\beta$ is an inner
perturbation of $\sigma\alpha\sigma^{-1}$ for some automorphism
$\sigma$ of $\mathfrak{B}$, where $\sigma\alpha\sigma^{-1}$  is the
action $t\mapsto\sigma\alpha_t\sigma^{-1}$.

If one of the above conditions are satisfied and $\alpha$ is
inner(approximately or asymptotically inner), then so is $\beta$.

{\bf Proof:} The first statement of the proposition can be obtained from [3,Corollary 1.3].
We just give the proof of the last statement of the proposition.

First we prove that if  $\beta$ is an inner perturbation of $\alpha$,
i.e. $\delta_\beta=\delta\alpha+i$ad$P$,  then it follows that if
$\alpha$ is inner(approximately or asymptotically inner), then so is
$\beta$.

(a) If $\alpha$ is inner,then so is $\beta$ by Corollary 2.5.

(b) If $\alpha$ is approximately inner, then there is a sequence $\{h_n\}$ in
$ \mathfrak{B}$ such that
$$\lim_{n\rightarrow\infty}\max_{|t|\leq 1} ||\alpha_t(A)-Ade^{th_n}(A)||=0.$$
$\alpha_{n,t}\triangleq Ade^{th_n}$, we construct $\alpha_n$-cocycle
$u_{n,t}$(resp. $u$)  for $\alpha_{n,t}$(resp. $\alpha$) such that
$\frac{d}{dt}u_{n,t}|_{t=0}=iP$ by Lemma 2.3. Then
since $\beta_t=Adu_t\circ\alpha_t$ and
$u_{n,t}\rightarrow u_t$,  we obtain that
$Adu_{n,t}\circ\alpha_{n,t}\rightarrow \beta_t$. Besides,
$Adu_{n,t}\circ\alpha_{n,t}$ is inner by Corollary 2.5. Then
$\beta$ is approximately inner.

(c) If $\alpha$ is asymptotically  inner, then
there is a continuous function $h$ of
$\mathbb{R}_+$
 into $\mathfrak{B}$ such that
 $$\lim_{s\rightarrow \infty}\max_{|t\leq 1|}\parallel\alpha_t(A)-Ade^{th(s)}(A)\parallel=0$$
 for any $A\in\mathfrak{B}$.

 $\alpha_{s,t}\triangleq Ade^{th(s)}$, we construct $\alpha_s$-cocycle
$u_{s,t}$(resp. $u$)  for $\alpha_{s,t}$(resp. $\alpha$) such that
$\frac{d}{dt}u_{s,t}|_{t=0}=iP$ by Lemma 2.3. Then
since $\beta_t=Adu_t\circ\alpha_t$ and
$u_{s,t}\rightarrow u_t$,  we obtain that
$Adu_{s,t}\circ\alpha_{s,t}\rightarrow \beta_t$. Besides,
$Adu_{s,t}\circ\alpha_{s,t}$ is inner and $Adu_{s,t}\circ\alpha_{s,t}=Ade^{t(h(s)+P)}$ by Corollary 2.5. Then
$\beta$ is approximately inner.

Next we shall prove that if $\beta$ is conjugate to $\alpha$, i.e. there is a
bounded automorphism $\sigma$ of $\mathfrak{B}$ such that
$\beta=\sigma\alpha\sigma^{-1}$. Then if $\alpha$ is
inner(approximately or asymptotically inner), then so is $\beta$.

(a') If $\alpha$ is inner, i.e. $\alpha_t(A)=e^{th}Ae^{-th}$. Then
$\beta_t(A)=\sigma(\alpha_t(\sigma(A))) =e^{t\sigma(h)}Ae^{-t\sigma(h)}$
i.e. $\beta$ is inner.

(b') If $\alpha$ is approximately inner, then there exits a sequence $\{h_n\}$ in
$ \mathfrak{A}$ such that
$$\lim_{n\longrightarrow\infty}\max_{|t|\leq 1} ||\alpha_t(A)-Ade^{th_n}(A)||=0.$$
Because $\beta$ is conjugate to $\alpha$, i.e. $\beta=\sigma\alpha\sigma^{-1}$,
$\sigma$ is the bounded automorphism  of $\mathfrak{A}$, so there exists an $M>0$ such
that $||\sigma||\leq M, ||\sigma^{-1}||\leq M$.

Therefore,
$$||\beta_t(A)-Ade^{t\sigma(h_n)}(A)||=
||\sigma^{-1}(\beta_t(\sigma(A)))-\sigma^{-1}(e^{t\sigma(h_n)}(\sigma(A)))||$$

i.e.$\alpha$ is approximately inner.

(c') If $\alpha$ is asymptotically inner, a similar argument shows
that $\beta$ is asymptotically inner.

Finally, if $\beta$ is cocycle-conjugate to $\alpha$, then there is a bounded
automorphism $\sigma$ of $\mathfrak{B}$ such that
$\beta$ is a cocycle perturbation of $\sigma\alpha\sigma^{-1}$. If $\alpha$
is inner(approximately or asymptotically inner), then so is $\sigma\alpha\sigma^{-1}$,
then so is $\beta$.

\vskip1cm
\section{Cocycle perturbations on nest algebra\\ and quasi-triangular algebra}

In this section $\mathfrak{H}$ denotes a Hilbert space,
 $\mathfrak{B}$ denotes a nest algebra (or a quasi-triangular algebra ) on $\mathfrak{H}$.
  $id$ denotes the identity automorphism on $\mathfrak{B}$ .

Let us mention briefly some well-known results on this subject.

{\bf Proposition 3.1$^{[11]}$:} If $\sigma$ is an automorphism on
$\mathfrak{B}$ and
$||\sigma-id||<1$,then there is an invertible element $T$ in
$\mathfrak{B}$ such
that $\sigma(A)=TAT^{-1}$ and $||T-I||<4||\alpha-id||$.

{\bf Proposition 3.2:} Let $\alpha$ be a flow
on $\mathfrak{B}$ and
$||\alpha_t-id||=O(t)$. Then there is a $P\in\mathfrak{B}$
such that $\alpha_t(A)=e^{tP}Ae^{-tP}$,
$\forall A\in\mathfrak{B}$.

{\bf Proof:}  When $||\alpha_t-id||=O(t)$,  there is an
$\varepsilon>0$ so that $|t|<\varepsilon$ implies
$||\alpha_t-id||<1$. By Proposition 3.1, there are invertible
elements $\{T_t\}$ in $\mathfrak{B}$ such that
$\alpha_t(A)=T_tAT_t^{-1}$ and $||T_t-I||<4||\alpha_t-id||$. By [1,
Proposition 3.1.1] $t^{-1}||\alpha_t-id||$ is bounded. So
$P_n=n(T_{\frac{1}{n}}-I)$ is bounded. So there is a weak$^*$
convergent sub-net $\{P_{n_r}\}$ with limit $P$ in $\mathfrak{B}$.
Let $\delta_\alpha$ denote the infinitesimal generator of $\alpha$,
one obtains
$$\delta_\alpha(A)=\lim_{r\rightarrow\infty}n_r(\alpha_{\frac{1}{n_r}}-id)AT_{\frac{1}{n_r}}=
\lim_{r\rightarrow\infty}
n_r(T_{\frac{1}{n_r}}A-AT_{\frac{1}{n_r}})=
\lim_{r\rightarrow\infty}P_{n_r}A-AP_{n_r}=PA-AP$$
for $\forall A\in\mathfrak{B}$.
Then by [1, Proposition 3.1.1]  $\alpha_t(A)=e^{tP}Ae^{-tP}$.

{\bf Proposition 3.3:} Let $\alpha, \beta$ be two flows on
$\mathfrak{B}$ and
$||\alpha_t-\beta_t||=O(t)$, then $\alpha$ is an inner perturbation of $\beta$.

{\bf Proof:} There exists an $M \geq 1$ and $\xi\geq$inf$_{t\neq
0}$($t^{-1}$log$||\alpha_t||$) such that
$$||\alpha_t||\leq Me^{\xi |t|}, ||\beta_t||\leq Me^{\xi |t|} $$
Since  $||\alpha_t\circ\beta_{-t}-id||\leq Me^{\xi |t|}
||\alpha_t-\beta_t||=O(t)$, we get that there is a
$P\in\mathfrak{B}$ such that
$\alpha_t(A)=e^{tP}\beta_t(A)e^{-tP}$, for
$A\in\mathfrak{B}$ by Proposition
3.2.
Finally $\delta_\alpha=\delta_\beta+$ad$i(-iP)$ is straightforward.
i.e. $\alpha$ is an inner perturbation of $\beta$.

{\bf Remark  :} By the definition of cocycle perturbation, it is obvious
that if $\beta$ is an cocycle perturbation of $\alpha$, then
$||\alpha_t-\beta_t||=O(t)$. Hence, cocycle perturbation and inner
perturbation are equivalent on $\mathfrak{B}$ by Lemma 2.3 and Proposition 3.3.
From the next theorem we obtain that
every flow on nest algebra or quasi-triangular algebra are uniformly
continuous.

{\bf Theorem 3.4:} Every flow $\mathfrak{B}$ are uniformly continuous.
And any two flows on $\mathfrak{B}$ are cocycle perturbation
to each other.

{\bf Proof:} In fact, if there is a flow $\alpha$ on
 $\mathfrak{B}$ which is not
uniformly continuous, let $\delta_\alpha$ be the infinitesimal generator of
$\alpha$.  Then there is a $v\in G(\mathfrak{B})$ such that
$v\notin D(\delta)$. Let
$v_t=v^{-1}\alpha_t(v)$, $\beta_t(A)=v_t\alpha_t(A)(v_t)^{-1}, t\in\mathbb{R}, A\in\mathfrak{B}$.
Then $\beta$ is an cocycle perturbation of $\alpha$. And $v_t$ is
not differentiable.

By Proposition 3.3  and Remark, $\beta$ is a inner perturbation of
$\alpha$ i.e. there is differentiable  $\alpha$-cocycle
$\{u_t\}\subseteq\mathfrak{B}$ such that
$\beta_t(A)=u_t\alpha_t(A){(u_t})^{-1}$.

But the commutant of  $\mathfrak{B}$ is trivial [11]. Hence it
follows that $v_t=\lambda_t u_t, t\in\mathbb{R}$,  where
$\lambda_t\in\mathbb{C}$.
 $u_t, v_t$ are $\alpha$-cocycle, it follow that $\lambda_{t+s}=\lambda_t\lambda_s$ and
$\lambda_t\rightarrow 1$, for $t\rightarrow 0$. Therefore
$\lambda_t$ is differentiable and $v_t$ is differentiable. This is
contrary to $v_t$ is
not differentiable.

Finally if $\alpha, \beta$ are flows on $\mathfrak{B}$, then
$\alpha, \beta$ are uniformly continuous. So $\alpha, \beta$
are inner perturbation to $\gamma_t\equiv id$. Because the inner perturbation defines
an equivalence relation on the flows, so $\alpha$ is an inner
perturbation of $\beta$, and $\alpha$ is a cocycle perturbation
of  $\beta$.

\nocite{liyk2,liyk/kua1}

\begin{thebibliography}{4}


\bibitem{B1979}
Bratteli, O. and Robinson, D. W.,  Operator Algebras and Quantum
Statistical Mechanics, I, Springer, 1979.


\bibitem{S1991}
Sakai, S., Operator algebras in dynamical systems, Cambridge Univ. Press, Cambridge 1991.



\bibitem{K2000}
 Kishimoto, A., Locally representable one-parameter automorphism
groups of AF algebras and KMS states, Rep. Math. Phys. 45 (2000), 333--356.



\bibitem{K2001}
 Kishimoto, A.,  UHF flows and the flip automorphism.
Rev. Math. Phys. 13 (2001), no. 9, 1163--1181

\bibitem{K2001}
Kishimoto, A.,  Examples of one-parameter automorphism
groups of UHF algebra. Commun. Math. Phys. 216 (2001), 395--428


\bibitem{K2002}
Kishimoto, A., Approximately inner flows on seperable $C^*$-algebras.
Rev. Math. Phs. 14 (2002), 1065--1094.


\bibitem{K2005}
Kishimoto, A., Approximate AF flows. J. Evol. Equ. 5 (2005), 153--184.


\bibitem{K2005}
Kishimoto, A., The one-cocycle property for shifts. Ergod. Th. And. Dynam. Sys. 25 (2005), 823--859.


\bibitem{K2006}
Kishimoto, A., Multiplier cocycles of a flow on a $C^*$-algebra. J. Funct. Anal. 235 (2006), 271--296.



\bibitem{K2006}
Kishimoto, A.,
 Lifting of an asymptotically inner flow for a separable $C\sp *$-algebra.
 Operator Algebras: The Abel Symposium 2004, 233--247,
Abel Symp., 1, Springer, Berlin, 2006.

 \bibitem{D1988}
 Davidson K. R., Nest algebras, Longman group UK limited, Essex, 1988.

\end{thebibliography}

\end{document}